\documentclass[reqno]{amsart}
\begin{document}
\title[]
{Periodic solutions of Abel differential equations}
\author[]
{M. A. M. Alwash}

\address{M. A. M. Alwash \hfill\break
Department of Mathematics,
West Los Angeles College,
9000 Overland Avenue, Los Angeles, CA 90230-3519, USA}
\email{alwashm@wlac.edu}

\subjclass[2000]{34C25, 34C07, 34C05, 37G15}
\keywords{Periodic solution; Abel differential equation; rigid system; 
\hfill\break\indent
limit cycle; Hopf bifurcation}

\begin{abstract}
For a class of polynomial non-autonomous differential equations of degree $n$, 
we use phase plane analysis to show that each equation in this class has 
$n$ periodic solutions. The result implies that certain rigid two-dimensional 
systems have at most one limit cycle which appears through multiple 
Hopf bifurcation.
\end{abstract}

\maketitle
\numberwithin{equation}{section}
\newtheorem{theorem}{Theorem}[section]
\newtheorem{lemma}[theorem]{Lemma}
\newtheorem{proposition}[theorem]{Proposition}
\newtheorem{remark}[theorem]{Remark}
\newtheorem{corollary}[theorem]{Corollary}

\section{Introduction}

We consider differential equations of the form 
\begin{equation} \label{e1.1}
\dot{z}:=\frac{dz}{dt}=z^{n}+P_{1}(t)z^{n-1}+\dots+
P_{n-1}(t)z+P_{n}(t) 
\end{equation} 
where $z$ is a complex-valued function
and $P_{i}$ are real-valued continuous functions. 
We denote by $z(t,c)$ the solution of \eqref{e1.1} satisfying $z(0,c)=c$.
Take a fixed real number $\omega$, we define the set $Q$ to be the set
of all complex numbers $c$ such that $z(t,c)$ is defined for all $t$ in the
interval $[0,\omega]$; the set $Q$ is an open set. On $Q$ we define the
displacement function $q$ by
$$
q(c)=z(\omega,c)-c.
$$
Zeros of $q$ identify initial points of solutions of \eqref{e1.1} which satisfy
the boundary conditions $z(0)=z(\omega)$. We describe such solutions as
{\em periodic} even when the functions $P_{i}$ are not themselves periodic.
However, if $P_{i}$ are $\omega$-periodic then these solutions are also
$\omega$-periodic. The main concern is to estimate the number of periodic 
solutions. This problem was suggested by C. Pugh as a version of Hilbert's 
sixteenth problem; it is listed by S. Smale as Problem 7 in \cite{s1}. 
Equations \eqref{e1.1} have been studied in detail by Lloyd in \cite{l2}, 
using the methods of complex analysis and topological dynamics. 

Note that $q$ is holomorphic on $Q$. The multiplicity of a periodic solution
$\varphi$ is that of $\varphi (0)$ as a zero of $q$. It is useful to work
with a complex dependent variable. The reason is that the number of zeros 
of a holomorphic function in a bounded region of the complex plane cannot 
be changed by small perturbations of the function. Hence, periodic solutions
cannot then be destroyed by small perturbations of the right-hand side of
the equation; periodic solutions can be created or destroyed only at infinity. 
Suppose that $\varphi$ is a periodic solution of multiplicity $k$. 
By applying Rouche's theorem to the function $q$, for any sufficiently small
perturbations of the equation, there are precisely $k$ periodic solutions in
a neighborhood of $\varphi$ (counting multiplicity). Upper bounds on the
number of periodic solutions of \eqref{e1.1} can be used as upper bounds on
the number of periodic solutions when $z$ is limited to be real-valued.
This is the reason that $P_{i}$ are not allowed to be complex-valued. 

When $n=3$, equation \eqref{e1.1} is known as the Abel 
differential equation. It was shown in \cite{l2} and \cite{p2} that 
Abel differential equation has exactly three
periodic solutions provided account is taken of multiplicity. We describe 
equation \eqref{e1.1} as of Abel form. 
For $n \geq 4$, Lins Neto \cite{l1} has given
examples which demonstrate that there is no upper bound, in terms of $n$ only, 
on the number of periodic solutions. However, there are upper bounds for 
certain classes of equations. It was shown in \cite{l2}, that each of the 
following equations has $n$ periodic solutions.
\[
\dot{z}=z^{n}+\alpha(t)z,
\]
\[
\dot{z}=z^{n}+\alpha(t)z^{\frac{n+1}{2}}+\beta(t)z, 
\]
where $n$ is odd in the second equation.
Ilyashenko, in \cite{i1}, gave the 
following upper bound for the number of periodic solutions
\[
8\exp{((3C+2) \exp{(\frac{3}{2}(2C+3)^{n})})}
\]
where, $C > 1$ is an upper bound for the absolute values of the coefficients  
$P_{i}(t)$. Although this estimate is non-realistic, it is the only known 
explicit estimate. In \cite{p1}, Panov considered the equation 
\[
\dot{z}=z^{n}+\alpha(t)z^{2}+\beta(t)z+\gamma(t) 
\]
and proved that the equation has at most three real periodic solutions if 
$n$ is odd. Quartic equations having at least ten periodic solutions were 
described in \cite{a1}; the coefficients were polynomial functions in $t$ 
of degree $3$.

In this paper, we consider the class of equations 
\begin{equation}
\dot{z}=z^{n}+\alpha(t)z^{n-1}+\beta(t)z^{n-2} \label{e1.2}
\end{equation}
as a generalized Abel differential equation. 
We show that if $\beta(t) \leq 0$, then this equation 
has at most two non-zero periodic solutions. We give conditions on 
$\alpha$ and $\beta$ that imply the equation has exactly two non-zero periodic 
solutions, one non-zero periodic solution, or no non-zero periodic solutions. 
Particular cases of this result, with $n=4$ and $n=5$, were given in 
\cite{a4} and \cite{a3}. In Section 2, we 
describe the phase portrait of \eqref{e1.1} and recall some results from 
\cite{l2}. In Section 3, we state and prove our main result. In the final 
Section, we use the result of Section 3 to show that a certain family of 
rigid two-dimensional systems has at most one limit cycle and this limit 
cycle appears through multiple Hopf bifurcation. 

\section{The Phase Portrait}

We identify equation \eqref{e1.1} with the $n-$tuples 
$(P_{1},P_{2},...,P_{n})$
and write $\mathcal{L}$ for the set of all equations of this form.
With the usual definitions of additions and scalar multiplications,
$\mathcal{L}$ is a linear space; it is a normed space if for 
$P=(P_{1},P_{2},...,P_{n})$ we define
 \[
 \|P\| = \max\{\max_{0\leq t \leq \omega}|P_{1}(t)|,
\max_{0\leq t \leq \omega}|P_{2}(t)|,\cdots, 
\max_{0\leq t \leq \omega}|P_{n}(t)|\}
\]
The displacement function $q$ is holomorphic on the open set $Q$.
Moreover, $q$ depends continuously on $P$ with the above norm on
$\mathcal{L}$ and the topology of uniform convergence on compact sets
on the set of holomorphic functions.

The positive real axis and the negative real axis are invariant.
Moreover, if $\varphi$ is a non-real solution which is periodic,
then so is $\bar{\varphi}$, its complex conjugate.

In \cite{l2}, it was shown that the phase portrait of \eqref{e1.1}
is as shown in Figure 1.  
We refer to \cite{l2} for the
details. There, the coefficients $P_{i}(t)$ were
$\omega-$periodic. It can be verified that the same methods are
applicable to the study of the number of solutions that satisfy
$z(0)=z(\omega)$ whether the coefficients are periodic or not.

\begin{figure}[ht]
\begin{center}
\begin{picture}(150,150)
\put(75,75){\circle{40}}
\put(90,75){\vector(1,0){10}}
\put(106,72){\vector(0,-1){10}}
\put(106,78){\vector(0,1){10}}
\put(120,100){\vector(-1,1){10}}
\put(100,115){\vector(1,-1){10}}
\put(100,30){\vector(1,1){10}}
\put(115,60){\vector(-1,-1){10}}
\put(95,95){\vector(-1,-1){10}}
\put(95,55){\vector(-1,1){10}}
\put(94,70){\line(1,0){35}}
\put(94,80){\line(1,0){35}}
\put(93,83){\line(1,1){35}}
\put(86,90){\line(1,1){35}}
\put(93,67){\line(1,-1){35}}
\put(86,59){\line(1,-1){35}}
\put(132,70){\makebox(0,0)[bl]{$G_{0}$}}
\put(130,130){\makebox(0,0)[bl]{$G_{1}$}}
\put(130,100){\makebox(0,0)[bl]{$H_{0}$}}
\put(90,125){\makebox(0,0)[bl]{$H_{1}$}}
\put(85,15){\makebox(0,0)[bl]{$H_{2n-4}$}}
\put(130,20){\makebox(0,0)[bl]{$G_{2n-3}$}}
\put(130,50){\makebox(0,0)[bl]{$H_{2n-3}$}}
\put(70,72){\makebox(0,0)[bl]{$D$}}
\thicklines
\end{picture}
\end{center}
\caption{Phase Portrait}
\end{figure}
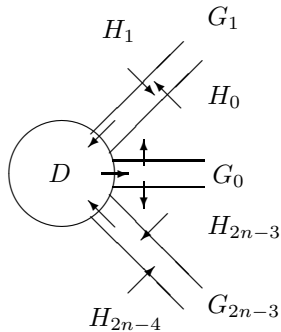

Note that the radius, $\rho$, of the disc $D$ is a sufficiently large number 
that depends only on $\|P\|$ and $\omega$. If $z=r e^{i\theta}$ then the sets
$G_{k},k=0,1,\dots,2n-3$, which are the "arms" in the figure, are
defined by
\[
G_{k}= \{z| r>\rho, \frac{k\pi}{n-1}-\frac{a}{r} < \theta 
< \frac{k\pi}{n-1}+\frac{a}{r}\}
\]
where $a=\max\{6,6\|P\|\}$. Between the arms are the sets
$H_{k},k=0,1,\dots,2n-3$, which are defined by
\[
H_{k}= \{z| r>\rho, \frac{k\pi}{n-1}+\frac{a}{r} \leq \theta 
\leq \frac{(k+1)\pi}{n-1}-\frac{a}{r}\}
\]
For even $k$, trajectories can enter $G_{k}$ only across $r=\rho$,
and for odd $k$, trajectories can leave $G_{k}$ only across $r=\rho$.
No solution can become infinite in $H_{k}$ as time either increases or
decreases. Every solution enters $D$. Solutions become unbounded
if and only if they remain in one of the arms $G_{k}$, tending to infinity 
as $t$ increases if $k$ is even and as $t$ decreases if $k$ is odd.

Let $q(P,c)=z_{P}(\omega,c)-c$, where $z_{P}(t,c)$ is the solution of
 $P \in \mathcal{L}$ satisfying $z_{P}(0,c)=c$. Suppose that $(P_{j})$ and
$(c_{j})$ are sequences in $\mathcal{L}$ and $\mathbb{C}$,
respectively, such that $q(P_{j},c_{j})=0$. If $P_{j} \to  P$ and
$c_{j} \to  c$ as $j \to  \infty$, then either $q(P,c)=0$, in this
case $z_{P}(t,c)$ is a periodic solution, or $z_{P}(t,c)$ is not
defined for the whole interval $0 \leq t \leq \omega$. In the
later case, we say that $z_{P}(t,c)$ is a {\em singular periodic
solution}. We also say that $P$ has a singular periodic solution
if $c_{j} \to  \infty$; in this case there are $\tau$ and $c$ such
that the solution $z_{P}$ with $z_{P}(\tau)=c$ becomes unbounded
in finite time as $t$ increases and as $t$ decreases. We summarize
the results of \cite{l2}.

\begin{proposition} \label{prop2.1}
\begin{itemize}
\item[(i)] Let $\mathcal{A}$ be the subset of $\mathcal{L}$ consisting of all
equations which have no singular periodic solutions.
The set $\mathcal{A}$ is open in $\mathcal{L}$. All equations in the same
components of $\mathcal{A}$ have the same number of periodic solutions.
\item[(ii)] The equation $\dot{z}=z^{n}$ has exactly $n$ periodic solutions.
\item[(iii)] For each $k$ and for each $t_{0}$ and $r \geq \rho$, 
there is a unique $\theta_{0}$ such
that the solution $z(t,c)$ with 
$z(t_{0},c)=c$ 
remains in $G_{k}$ for $t \geq t_{0}$ ($k$ is even) or 
$t \leq t_{0}$ ($k$ is odd) if and only if $c=re^{i \theta_{0}}$. 
\end{itemize}
\end{proposition}

\section{Main Result}

Assume that $P_{n}(t) \equiv 0$. We call the solution $z=0$ a {\em center} 
if $z(t,c)$ is periodic
for all $c$ in a neighborhood of $0$. If the term $z^{n}$ in \eqref{e1.1} is 
replaced by $P_{0}(t)z^{n}$, then there are equations with a center when 
$P_{0}$ has zeros. For cubic equations, this is
related to the classical center problem of polynomial
two-dimensional systems; we refer to \cite{a5} for details.
However, when $P_{0}$ has no zeros then $z=0$ is never a center.
Particular cases of this result were given in \cite{a3} and \cite{a4} for 
$n=4$ and $n=5$, respectively. We give a brief proof, for the sake of 
completeness.

\begin{theorem} \label{thm3.1}
The solution $z=0$ is isolated as a periodic solution of \eqref{e1.1} with 
$P_{n}(t) \equiv 0$.
\end{theorem}

\begin{proof}
Suppose, if possible, that there is a open set
$U \subset \mathbb{C}$ containing the origin such that all solutions
starting in $U$ are periodic. Then $q \equiv 0$ in the component of
its domain of definition containing the origin. But the real zeros
of $q$ are contained in the disc $D$. Thus
\[
\inf\{c \in \mathbb{R}: c>0, z(t,c)\mbox{ is not defined for }
0 \leq t \leq \omega\} < \infty
\]
It follows that there is a real singular periodic solution;
but a positive real periodic solution which tends to infinity can do
so only as $t$ increases. This is a contradiction, and the result follows.
\end{proof}

Now, we give the result about the number of periodic solutions.

\begin{theorem} \label{thm3.2}
 Suppose that $\beta(t) \leq 0$. Equation \eqref{e1.2} has exactly $n$ periodic 
solutions.
\end{theorem}

\begin{proof} With $z=re^{i\theta}$, we have
\[
\dot{\theta}=r^{n-1}\sin{(n-1)\theta}+r^{n-2}\alpha(t) \sin{(n-2)\theta}+
r^{n-3}\beta(t) \sin{(n-3)\theta}.
\]
If $|c| > \rho$ and is real then the real solution $z(t,c)$ remains outside
the disk $D$ either when $t$ increases or when $t$ decreases, 
and will become infinite.
Solutions that enter $G_{0}$ or $G_{n-1}$ will leave $G_{0}$ or $G_{n-1}$,
except the solution that enters at $\rho e^{i\theta_{0}}$ described in 
part (iii) of Proposition \ref{prop2.1}; 
this solution is real because any solution which is
once real is always real. Therefore, the unique solution that becomes
infinite is a real solution if $k=0$ or $k=n-1$. On the other hand, 
no real solution is unbounded as
$t$ increases and decreases. Hence, no singular periodic solution
enters $G_{0}$ or $G_{n-1}$ because singular periodic solutions are
unbounded both as $t$ increases and decreases.
Thus, a singular periodic solution enters $D$ from a $G_{k}$ with odd 
$k \neq n-1$ 
and leaves $D$ to a $G_{j}$ with even $j \neq 0, n-1$. 

Now, if $k$ is odd and $ 1 \leq k \leq n-3$, let 
$\theta_{1}=\frac{k\pi}{n-2}$, $\theta_{2}=\frac{(k-1)\pi}{n-2}$, and consider
\begin{gather*}
\dot{\theta}(\theta_{1})=r^{n-1}\sin{((n-1)\theta_{1})}+
r^{n-2}\alpha(t) \sin{((n-2)\theta_{1})}+r^{n-3}\beta(t) \sin{((n-3)\theta_{1})}
,\\
\dot{\theta}(\theta_{2})=r^{n-1}\sin{((n-1)\theta_{2})}+
r^{n-2}\alpha(t) \sin{((n-2)\theta_{2})}+r^{n-3}\beta(t) \sin{((n-3)\theta_{2})}
\end{gather*}

But, 
\[
\sin{((n-1)\theta_{1})}=\sin{(k\pi+\frac{k\pi}{n-2})} < 0, 
\]
\[
\sin{((n-3)\theta_{1})}=\sin{(k\pi-\frac{k\pi}{n-2})} > 0,
\]
\[
\sin{((n-1)\theta_{2})}=\sin{((k-1)\pi+\frac{(k-1)\pi}{n-2})}  > 0, 
\]
\[
\sin{((n-3)\theta_{2})}=\sin{((k-1)\pi-\frac{(k-1)\pi}{n-2})} < 0,
\]
\[\sin{((n-2)\theta_{1})}=\sin{((n-2)\theta_{2})}=0.
\]
Under the above hypotheses, $\dot{\theta}(\theta_{1}) < 0$ and
$\dot{\theta}(\theta_{2}) > 0$. Since $\rho$ is a sufficiently 
large number, we assume that $\rho > \frac{a(n-2)(n-1)}{\pi}$. This condition 
on $\rho$ guarantees that the arc of intersection of $G_{k}$ with $D$ lies 
inside the sector 
\[
\{r<\rho, \theta_{2} < \theta < \theta_{1} \}.
\]
Precisely, the condition implies that 
\[
\frac{k\pi}{n-2} > \frac{k\pi}{n-1}+\frac{a}{\rho},
\frac{(k-1)\pi}{n-2} < \frac{k\pi}{n-1}-\frac{a}{\rho}.
\]
Hence, solutions do not leave the sector (see Figure 2).
Therefore, no singular periodic
solution can enter $D$ from $G_{k}$ and leave $D$ to $G_{k-1}$ or $G_{k+1}$.
Since the phase portrait is symmetric about the $x-$axis, if $k$ is odd 
and $n \leq k \leq 2n-3$, no singular
periodic solution can enter $D$ from $G_{k}$ and leaves $D$
to $G_{k-1}$ or $G_{k+1}$. 
It follows that the equation does not have a singular
periodic solution.

Now, consider the class of equations
\[
\dot{z}=z^{n}+s \alpha(t)z^{n-1}+s \beta(t)z^{n-2},
\]
with $0 \leq s \leq 1$.
Since, $s \beta(t) \leq 0$, 
any equation in this family does not have singular
periodic solutions. The equation $\dot{z}=z^{n}$ belongs to this
family and has $n$ periodic solutions.
By part (i) of Proposition \ref{prop2.1}, each of these equations has 
$n$ periodic solutions.
\end{proof}

\begin{figure}[ht]
\begin{center}
\begin{picture}(150,150)
\put(75,75){\circle{40}}
\put(102,75){\vector(-1,0){10}}
\put(115,85){\vector(0,-1){10}}
\put(110,65){\vector(0,1){10}}
\put(82,90){\vector(1,1){10}}
\put(84,58){\vector(1,-1){10}}
\put(80,86){\vector(1,-1){10}}
\put(88,65){\vector(-1,2){5}}
\put(94,78){\line(1,0){35}}
\put(94,72){\line(1,0){35}}
\put(86,90){\line(1,1){35}}
\put(80,94){\line(1,1){35}}
\put(88,59){\line(1,-1){35}}
\put(82,56){\line(1,-1){35}}
\put(75,75){\line(2,1){35}}
\put(75,75){\line(2,-1){35}}
\put(132,70){\makebox(0,0)[bl]{$G_{k}$}}
\put(122,130){\makebox(0,0)[bl]{$G_{k+1}$}}
\put(125,12){\makebox(0,0)[bl]{$G_{k-1}$}}
\put(135,90){\makebox(0,0)[b1]{$\theta=\frac{k\pi}{n-2}$}}
\put(135,47){\makebox(0,0)[b1]{$\theta=\frac{(k-1)\pi}{n-2}$}}
\thicklines
\end{picture}
\end{center}
\caption{Sector around $G_{k}$, odd $k$}
\end{figure}
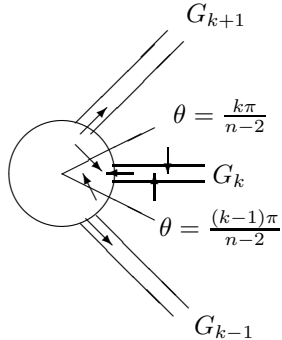

\begin{lemma} \label{lemma3.3}
Let $A=\int_{0}^{\omega} \alpha(t) dt$ and 
$B=\int_{0}^{\omega} \beta(t) dt$. The solution $z=0$ of \eqref{e1.2} 
has multiplicity 
\begin{itemize}
\item[(i)] $n-2$ if $B \neq 0$ 
\item[(ii)] $n-1$ if $B=0$ and $A \neq 0$.
\item[(iii)] $n$ if $A=B=0$.
\end{itemize}
\end{lemma}
\begin{proof}
We write 
\[
z(t,c)=\sum_{k=1}^{\infty} a_{k}(t) c^{k}\]
and substitute directly into the 
equation \eqref{e1.2}. This gives a recursive set of linear differential 
equations for the $a_{k}(t)$ with initial conditions $a_{1}(0)=1$ and 
$a_{k}(0)=0$ if $k >1$. The multiplicity is $K$ if $a_{1}(\omega)=1$, 
$a_{k}(\omega)=0$ for $2 \leq k \leq K-1$ and $a_{K}(\omega) \neq 0$. 
Direct computations give
\[
a_{1}(t) \equiv 1;  a_{k}(t) \equiv 0, 2 \leq k \leq n-3,
\]
\[
\dot{a}_{n-2}=\beta; \dot{a}_{n-1}=\alpha; \dot{a}_{n}=1.
\]
Solving the last three  equations gives $a_{n-2}(\omega)=B$, 
$a_{n-1}(\omega)=A$, and $a_{n}(\omega)=\omega$. The result follows.
\end{proof}
\begin{corollary} \label{coro3.4}
Consider equation \eqref{e1.2} with $\beta(t) \leq 0$. The equation has 
\begin{itemize}
\item[(i)] two non-zero periodic solutions if $B < 0$; at most one is 
positive and at most one is negative. If $B$ is small and negative, then 
there are two non-zero real periodic solutions,
\item[(ii)] one non-zero periodic solution if $B=0$ and $A \neq 0$; this 
solution is a real solution. It is positive if $A<0$ and it is negative 
if $A>0$, 
\item[(iii)] no non-zero periodic solutions if $A=B=0$.
\end{itemize}
\end{corollary}
\begin{proof}
Consider the case $B<0$. The multiplicity of the origin is $n-2$. There are 
two non-zero periodic solutions. In a neighborhood $c=0$, 
$q(c)=Bc^{n-2}+O(c^{n-1})$. 
If $n$ is odd then $q(c)$ has the sign of $c$ if $|c|$ is large. 
Therefore, $q$ has a positive solution and a negative solution if it is defined 
for large $|c|$. If $n$ is even then $q(c) > 0$ if $|c|$ is large. Again, 
$q$ has a positive solution and a negative solution if it is defined for 
large $|c|$. 

In the case $B=0$ and $A \neq 0$, the multiplicity of the origin is $n-1$. 
There is only one non-zero solution which is a real solution because 
complex solutions occur in conjugates pairs. The argument used in the 
first case implies that this solution is positive if $A<0$ and is negative 
if $A>0$.

If we start with $A=B=0$, then the multiplicity of $z=0$ is $n$. Perturb the 
equation, so that $B=0$ and $A \neq 0$. Since the total number of periodic 
solutions is unchanged by small perturbations, a real solution will bifurcate 
out of the origin. Now, we make a second perturbation so that $B$ is 
negative and such that $|B|$ is small compared to $|A|$. A second real 
nonzero solution will bifurcate out of the origin. If the equation is 
perturbed such that $A=0$ but $B$ is negative and small, then the stability 
will be reversed and two non-zero real periodic solutions bifurcate out of 
the origin; one is positive and one is negative.
\end{proof}

\section{RIGID SYSTEMS}

Consider the system
\begin{equation}
\begin{gathered}
\dot{x}=\lambda x - y +x(R_{n-1}(x,y)+R_{n-2}(x,y)+\dots+R_{1}(x,y))\\
\dot{y}=x+ \lambda y +y(R_{n-1}(x,y)+R_{n-2}(x,y)+\dots+R_{1}(x,y)),
\end{gathered}\label{e4.1}
\end{equation}
where $R_{i}$ is a homogeneous polynomial of degree $i$. The system
in polar coordinates becomes
\begin{gather*}
\dot{r}=r^{n}R_{n-1}(\cos{\theta},\sin{\theta})+
r^{n-1}R_{n-2}(\cos{\theta},\sin{\theta})+\dots+
r^{2}R_{1}(\cos{\theta},\sin{\theta})+\lambda r\\
\dot{\theta}=1.
\end{gather*}
This system is called a rigid system because the derivative of the angular 
variable is constant. 
It was shown in \cite{a2} that if the origin is a center then 
$\int_{0}^{2\pi} R_{k}(\cos{\theta},\sin{\theta}) d \theta =0$ for all 
$ 1 \leq k \leq n-1$. In fact, these definite integrals are the first 
focal values of the system. It is clear
that the origin is the only critical point and if it is a center then
it is a uniformly isochronous center; see \cite{c1}. 
Limit cycles of \eqref{e4.1} correspond to positive $2\pi-$periodic solutions of
\[
\frac{dr}{d\theta}=R_{n-1} r^{n} + R_{n-2} r^{n-1} + \dots
+ R_{1} r^{2} + \lambda r
\]
In \cite{g1}, the family of rigid systems 
\[
\dot{x}=\lambda x - y +xR_{n-1}(x,y), 
\]
\[
\dot{y}=x+ \lambda y +yR_{n-1}(x,y)
\]
was considered. It was shown that if 
$ \lambda \int_{0}^{2\pi} R_{n-1}(\cos{\theta},\sin{\theta}) d\theta <0$, then 
there is at most one limit cycle. Cubic rigid systems were considered in 
\cite{g2}. It was shown that there are at most two limit cycles. 

Now, we consider the rigid system
\begin{equation}
\begin{gathered}
\dot{x}= - y +x(R_{n-1}(x,y)+R_{n-2}(x,y)+R_{n-3}(x,y))\\
\dot{y}=x+ y(R_{n-1}(x,y)+R_{n-2}(x,y)+R_{n-3}(x,y)).
\end{gathered}\label{e4.2}
\end{equation}
In polar coordinates, the system becomes
\[
\frac{dr}{d\theta}=R_{n-1} r^{n} + R_{n-2} r^{n-1} + R_{n-3} r^{n-2}.
\]
If a function $R_{i}$ does not change sign, then it is necessary
to assume that $i$ is even. We assume that $n$ is odd. Hence, $R_{n-2}$ is 
a homogeneous polynomial in $\sin{\theta}$ and $\cos{\theta}$ of odd degree; 
Therefore, $\int_{0}^{2\pi} R_{n-2}(\cos{\theta},\sin{\theta}) d\theta=0$. 
On the other hand, real periodic solutions occur in pairs. 
If $\varphi(\theta)$ is a real periodic solution then so is 
$-\varphi(\theta+\pi)$. Let 
$B= \int_{0}^{2\pi} R_{n-3}(\cos{\theta},\sin{\theta}) d\theta$. The 
Liapunov quantities of the system are $B$ and $1$. If $B=0$ then the 
origin is unstable. By perturbing the coefficients of $R_{n-3}$ such that 
$B < 0$ and is small enough, the origin will be unstable and a 
small-amplitude limit cycle appears through multiple Hopf bifurcation. 
The following result follows directly from these remarks and 
Theorem \ref{thm3.2}.

\begin{theorem} \label{thm4.1}
Suppose that $n$ is odd, $R_{n-1} \equiv 1$, and $R_{n-3} \leq 0$. 
Let $B= \int_{0}^{2\pi} R_{n-3}(\cos{\theta},\sin{\theta}) d\theta$.
\begin{itemize}
\item[(i)] The origin is not a center for system \eqref{e4.2}.
\item[(ii)] If $B<0$ then system \eqref{e4.2} has at most one limit cycle. If 
$B$ is small enough then the system has a unique limit cycle; this 
limit cycle is unstable. 
\item[(iii)] If $B=0$ then system \eqref{e4.2} does not have 
a limit cycle.
\end{itemize}
\end{theorem}

\begin{remark} \label{rmk5.3} \rm
 If the leading coefficient $R_{n-1}$ does not vanish anywhere then
the transformation of the independent variable
\[
\theta \mapsto \exp(\int_{0}^{\theta} R_{n-1}(\cos{u},\sin{u}) du)
\]
reduces the polar equation into a similar equation but with a leading
coefficient equals one.
\end{remark}

\subsection*{Acknowledgments}
I am very grateful to the Department of Mathematics,
University of California, Los Angeles for the hospitality.

\end{document}